%% file: MainTIR_AC.tex
\pdfoutput=1
\documentclass{birkjour}
 \theoremstyle{definition}
 
 \theoremstyle{remark}

 \numberwithin{equation}{section}

 \newtheorem{theorem}{Theorem}[section]
 \newtheorem{corollary}[theorem]{Corollary}
 \newtheorem{lemma}[theorem]{Lemma}
 
 \theoremstyle{definition}
 \newtheorem{definition}[theorem]{Definition}
 \theoremstyle{remark}
 \newtheorem{remark}[theorem]{Remark}
 
 \numberwithin{equation}{section}

\begin{document}

%-------------------------------------------------------------------------
% editorial commands: to be inserted by the editorial office
%
%\firstpage{1} \volume{228} \Copyrightyear{2004} \DOI{003-0001}
%\seriesextraline{This is the Concrete Title of this Book\br H.E. R and S.T.C. W, Eds.}
%\firstpage{1}

\title[Minimal non-extensible precolorings and implicit-relations]
 {Minimal non-extensible precolorings and \\ implicit-relations}

%----------Author 1
\author[Jos\'e Antonio Mart\'in H.]{Jos\'e Antonio Mart\'in H.}

\address{Faculty of Computer Science, Complutense University of Madrid, Spain.}
\email{jamartinh@fdi.ucm.es}

%\thanks{This work was completed with the support of our \TeX-pert.}



\date{\today}
%\dedicatory{To Julio Subocz in memoriam}
%%% ----------------------------------------------------------------------

\begin{abstract}
\input{AbstractTIR.tex}

\end{abstract}

%%% ----------------------------------------------------------------------
\maketitle
%%% ----------------------------------------------------------------------
%\tableofcontents

\input{BodyTIR.tex}

%\subsection*{Acknowledgment}
%Many thanks to our \TeX-pert for developing this class file.

\bibliographystyle{plain}
\bibliography{gt}

% ------------------------------------------------------------------------
\end{document}

%% file: AbstractTIR.tex
In this paper I study a variant of the general vertex coloring problem called precoloring. Specifically, I study graph precolorings, by developing new theory, for characterizing the minimal non-extensible precolorings. It is interesting per se that, for graphs of arbitrarily large chromatic number, the minimal number of colored vertices, in a non-extensible precoloring, remains constant; only two vertices $u,v$ suffice. Here, the relation between such $u,v$ is called an implicit-relation, distinguishing two cases: (i) implicit-edges where $u,v$ are precolored with the same color and (ii) implicit-identities where $u,v$ are precolored distinct.

%% file: BodyTIR.tex
\section{Introduction}

%A theory about the implication structure in graph coloring is presented. This theory is based on proposed concept of implicit-relations: implicit-edges and implicit-identities. This study may help to unveil a wide range of ``hidden" phenomena in a range of scientific disciplines such as physics, chemistry, and of course, discrete mathematics and computer science.

In this work, a formal definition of a theory of graph-chromatic implicit-relations is presented. This theory can be described a specialization of a variant of the general vertex coloring problem called precoloring, specifically, for characterizing minimal non-extensible precolorings. For a good reference on graph coloring (a.k.a. chromatic graph theory) the reader can see the introductory book of Chartrand and Zhang~\cite{Chartrand2008}, and also the book of Jensen and Toft~\cite{JT95} with interesting open problems. Chartrand and Zhang~\cite[pp. 240]{Chartrand2008} open the precoloring section affirming that:
\begin{quote}
...What we are primarily interested in, however, is whether a $k$-precoloring of $G$, where $k \geq \chi(G)$,
can be extended to a $k$-coloring of $G$...
\end{quote}
\noindent
However, in this paper, I show that we should pay attention, not only to $(k \geq \chi(G))$-precolorings of $G$, but, just to the 2-precolorings of $G$, since this are the minimal non-extensible precolorings, without regard of $\chi(G)$. This work is also relevant in connection with all theorems and problems exposed by Chartrand and Zhang, in the section about ``Precoloring Extensions of Graphs''~\cite[Chap.9, Sec.3]{Chartrand2008}.

It is easy to show that for an arbitrarily large chromatic number $k$, the minimal number of colored vertices, in a non-extensible precoloring, for some $k$-chromatic graph $G$, remains constant since only two vertices $u,v$ suffice (however it remains open and, up to my knowledge, there is no study about the minimal number of vertices of a non-extensible precoloring for some given particular graph).

We call the relation between such $u,v$ an implicit-relation distinguishing two cases:
\begin{description}
   \item[(i)]  implicit-edges, where $u,v$ has the same color, and this is a minimal non-extensible precoloring.
   \item[(ii)] implicit-identities, where $u,v$ has distinct color, and this is a minimal non-extensible precoloring.
\end{description}

Additionally, the notion of implicit-edge can be derived from the observation that, when finding a $k$-coloring of some $k$-chromatic graph $G$, there may be some independent set $S \subset V(G)$ such that there is no $k$-coloring of $G$ where all the vertices in $S$ receive the same color. Thus, when $|S|=2$, we say that the 2-element subset $S=\{u,v\}$ is an implicit-edge, based on the intuition that the $k$-colorings of $G$ are the same as if $uv \in E(G)$. Conversely, implicit-identities are, in most cases, the opposite-equivalent of implicit-edges. It is a relation in which vertices must receive always the same color being intuitively the same as if $G=G-uv/u,v$ (the vertex identification of $u$ and $v$).

Hence, given a $k$-chromatic graph $G$:
\begin{description}
   \item[(iii)] there is an implicit-edge when $u,v$ receive \emph{different} colors in every $k$-coloring of $G$.
   \item[(iv)]  there is an implicit-identity when $u,v$ receive \emph{the same} color in every $k$-coloring of $G$.
\end{description}

Implicit-relations have direct application, for instance, in the study of two very important open problems in graph theory, namely: the double-critical graphs problem and the Hadwiger's conjecture.

In short, a way to express the Hadwiger's conjecture is the following: every $k$-chromatic graph can be contracted (by successive edge contractions or edge and vertex deletions) into a complete graph on $k$ vertices ($K_k$). An equivalent formulation is to show that there are no contraction-critical graphs (i.e. a graph where each minor has lower chromatic number) different from $K_k$. The relation with implicit-edges is clear since contraction critical graphs cannot have minors with implicit-edges in the edge set. Hence a deep study of this concept is a way to approach the conjecture.

In this paper, in order to study and characterize such concepts, a series of theorems are presented.

\section{Preliminary definitions and basic terminology}
\label{sec:termninology}%Graph Coloring

Unless we state it otherwise, all graphs in this work are connected and simple (finite, and have
no loops or parallel edges).

Partitioning the set of vertices $V(G)$ of a graph $G$ into separate classes, in such a way that no two adjacent vertices are
grouped into the same class, is called the vertex graph coloring problem. In order to distinguish such classes, a set of colors C
is used, and the division into these \emph{(color) classes} is given by a proper-coloring (we will use here just the single term coloring) $c : V(G)\rightarrow \{1...k \}$, where $c(u) \neq c(v)$ for all $uv's$ belonging to the set of edges $E(G)$ of $G$. A \emph{k-coloring} of $G$ is a coloring that uses exactly $k$ colors. The \emph{Chromatic number} of a graph $\chi(G)$ is the minimum number such that there is a $\chi(G)$-coloring of a graph $G$. Thus, if $\chi(G) \leq k$ then one says that $G$ is \emph{k-colorable} (i.e. $G$ can be colored with $k$ different colors) and if $\chi(G)= k$ then one says that $G$ is \emph{k-chromatic}.

A \emph{precoloring} of a graph $G$, is a coloring $p: W \rightarrow \{1...k\}$ of $W \subset V(G)$ such that $p(u)\neq p(v)$ if $u,v\in W$ and $uv\in E(G)$. A precoloring $p$ of $G$ can be \emph{extended} to a coloring of all the vertices of $G$ when there is at least one coloring $c : V(G)\rightarrow \{1...k\}$ of $G$ such that $c(u)=p(u)$ for all $u \in W$, otherwise it is said that $p$ is non-extensible. A $k$-\emph{precoloring} $p$ is a precoloring of $G$ such that $p$ uses only $k$ colors.

An \emph{independent set} (also called \emph{stable set}) $I=\{u,v,w,...\}$  is a set of vertices of a graph $G$ such that there are no edges between any two vertices in $I$, i.e, if $\{u,v\}\in I$ then $uv\notin E(G)$. The set $I(G)$ will denote the set of all independent sets of graph $G$.

The set off all adjacent vertices to a vertex $u\in V(G)$ is called its \emph{neighborhood} and is denoted by $N_G(u)$ (when it is clear to which graph we are referring to, we will use simply $N(u)$, i.e. omitting the graph). The \emph{closed neighborhood} of a vertex $u$, denoted by $N[u]$, includes also the vertex $u$, i.e. $N(u)\cup \{u\}$.

The \emph{degree} of a vertex $u$, $deg(u)$, is equal to the cardinality of its neighborhood $deg(u)=|N(u)|$. A \emph{complete vertex} is any $u \in V(G)$ such that $N[u] = V(G)$ and a graph is called a \emph{complete graph} if every vertex is a complete vertex.

A \emph{path} $P_n$, of order $n$ and length $\ell$, is a sequence of $n$ joined vertices such that the travel from the \emph{start} vertex to the \emph{end} vertex passes trough $\ell$ edges. A \emph{simple path} has no repeated vertices.

Vertex deletions and additions are denoted as $G-u$ and $G+u$ respectively, for a graph $G$ and a vertex $u$. Edge deletions and additions are denoted as $G-uv$ and $G+uv$, or simply $G-e$ or $G+e$ respectively, for a graph $G$ and edge $e=uv$.

A \emph{vertex identification}, denoted by $G-uv/u,v$, is the process of replacing two \emph{non-adjacent} vertices $u,v$ of a graph $G$, i.e $uv\notin E(G)$, by a new vertex $w$ such that $N(w) = N(u)\cup N(v)$.

An \emph{edge contraction}, denoted by $G/uv$ or $G/e$, is the process of replacing two \emph{adjacent} vertices $u,v$ of G, i.e $uv \in E(G)$, by a new vertex $w$ such that $N(w) = N(u)\cup N(v)$.

A \emph{vertex contraction}, denoted by $G/u,v$, is the process of replacing two vertices $u,v$ of a graph $G$, by a new vertex $w$ such that $N(w) =N(u) \cup N(v)$. Hence \emph{vertex contractions} include both cases: vertex identifications and edge contractions.

A graph $H$ is called a \emph{minor} of the graph $G$, denoted as $H\prec G$, if $H$ is isomorphic to a graph that can be obtained from a subgraph of $G$ by zero or more edge deletions, edge contractions or vertex deletions on a subgraph of $G$. In particular, $G$ is minor of itself.

An element $x$ of a graph $G$ is called \emph{critical} if $\chi(G-x)<\chi(G)$. If all the vertices of a graph $G$ are critical we say that $G$ is \emph{vertex-critical} and if every element (vertex or edge) of $G$ is critical we say that $G$ is a \emph{critical graph} and more specifically if $\chi(G)=k$ we say that $G$ is $k$-critical.

An \emph{edge-subdivision} of an edge $e$ of a graph $G$ is the subdivision of some $e\in E(G)$ with endpoints $\{u,v\}$ resulting in a graph containing one new vertex $w$, and with an edge set replacing $e$ by two new edges $uw$ and $wv$ i.e., $H=G-e+w+uw+wv$.

Given a (pre)coloring $c : V(G)\rightarrow \{1...k\}$ of a graph $G$, a \emph{2-color-chain} $\Omega_{uv}$ (a.k.a \emph{Kempe's chain}) is a simple maximal bipartite connected subgraph $B \subset G$ such that $\{u,v\}\in V(B)$ and every vertex in $B$ has either the color $c(u)$ or $c(v)$. \emph{Flipping} a chain $\Omega_{uv}$ is the process of swapping the color of vertices with color $c(u)$ to $c(v)$ and, respectively, vertices with color $c(v)$ to $c(u)$.

\section{Implicit relations} \label{sec:main}

The chromatic implicit relations are mainly defined by their two most basics
concepts, implicit-edges and implicit-identities.

The notion of implicit-edge is based on the observation that, when finding a $k$-coloring of some $k$-chromatic graph $G$, there may be some independent set $S \in I(G)$ such that there is no $k$-coloring of $G$ where all the vertices in $S$ receive the same color. Thus, when $|S|=2$ we say that the 2-elements set $S=\{u,v\}$ is an implicit-edge. Respectively, implicit-identities are, in most cases, the opposite-equivalent of implicit-edges. While implicit-edges are defined as a relation in which vertices receive always different color implicit-identities are the opposite relation, that is, the vertices receive always the same color.

% -- definitions, trivial cases and graphs.
\begin{definition} \label{IE} Given a $k$-chromatic graph $G$, we say that $\{u,v\} \in V(G)$ is an \emph{implicit-edge} \textbf{iff} the set of all $k$-colorings of $G-uv$ where $u$ and $v$ receive the same color is the empty set:
\begin{equation}
\left\{c \in \Phi(G-uv) \; | \;  c(u)=c(v)\right\} = \emptyset,
\end{equation}
where $\Phi(G-uv)$ is the set of all $k$-colorings of $G-uv$ and $c(u),c(v)$ are the colors of vertices $u,v$ respectively, as assigned by a particular coloring ($c$) of $G$.
\end{definition}

\begin{remark}
Also, we must note that an implicit-edge could belong to the set of edges of $G$ or not, i.e. $\{u,v\}$ is also an implicit-edge in the graph $G-ij$. This can also be viewed as a precoloring that can't be extended to all $G$ (i.e. a non-extensible precoloring) by taking any $p: W \rightarrow C$ of $W \subset V(G)$ such that $p(u)=p(v)$.
\end{remark}

\begin{definition} \label{IE} Given a $k$-chromatic graph $G$, we say that $\{u,v\} \in V(G)$ is an \emph{implicit-identity} \textbf{iff} the set of all $k$-colorings of $G-uv$ where $u$ and $v$ receive different color is the empty set:
\begin{equation}
\{c \in \Phi(G-uv) \; | \;  c(u)\neq c(v)\} = \emptyset,
\end{equation}
where $\Phi(G-uv)$ is the set of all $k$-colorings of $G-uv$ and $c(u),c(v)$ are the colors of vertices $u,v$ respectively, as assigned by a particular coloring ($c$) of $G$.
\end{definition}
\begin{remark}
Note that contrary to implicit-edges, if we add the edge $e=uv$ to $G$ the resulting graph $G+e$ will not be $k$-colorable. This can also be viewed as a 2-precoloring that can't be extended to all $G$ by taking any $p: W \rightarrow \{1..k\}$ of $W \subset V(G)$ such that $p(u)\neq p(v)$.
\end{remark}

The trivial case of implicit-edges arise in bipartite (2-chromatic) graphs.
\begin{lemma}
Given a bipartite graph $G$, $\{u,v\}\in V(G)$ is an implicit-edge \textbf{iff} there is an odd-path $P_{2n}$ starting at $u$ and ending at $v$. (without provided proof)
\end{lemma}

Let us take the path $P_4$, as an illustrative example:
\begin{equation}
P_4=\; \mbox{$u\circ$---$\bullet$---$\circ$---$\bullet v$},
\end{equation}
where $\{u,v\}$ is an implicit-edge of $P_4$. If we identify vertices $\{u,v\}$ we get the $K_3$ graph. Also if we add the edge $uv$ to $P_4$ we get a square $C_4$ and then all edges become implicit-edges.

\begin{figure}[h]
\begin{centering}
  \begin{tabular}{ccc}
  \includegraphics[scale=1]{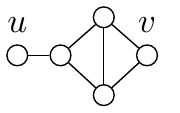} & \includegraphics[scale=1]{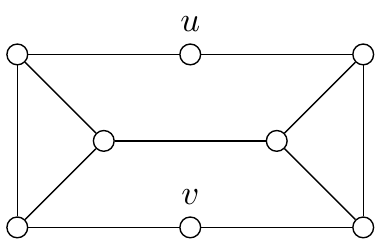} & \includegraphics[scale=1]{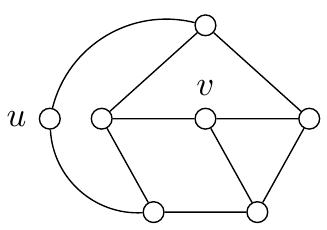} \\
  (A) & (B) & (C)
 \end{tabular}
  \caption{Implicit edges ($u,v$) in 3-chromatic planar graphs.}
  \label{impl1}
\end{centering}
\end{figure}

More complex examples are shown in Fig.~\ref{impl1}. for 3-chromatic planar graphs. Vertices $\{u,v\}$ forms an implicit edge since a 1-precoloring $p$ of $\{u,v\}$, i.e.  $p(u)=p(v)$, can't be extended to a 3-coloring of $G$.

Implicit-identities in bipartite graphs are defined very easy
\begin{lemma}
Given a bipartite graph $G$, $\{u,v\}\in V(G)$ is an implicit-identity \textbf{iff} there is an even-path $P_{2n+1}$ starting at $u$ and ending at $v$. (without provided proof)
\end{lemma}
Let us take the path $P_5$, as an illustrative example:
\begin{equation}
P_5=\; \mbox{$u\circ$---$\bullet$---$\circ$---$\bullet$---$\circ v$},
\end{equation}
where $\{u,v\}$ is an implicit-identity of $P_5$. If we add the edge $uv$ then we obtain a 3-chromatic graph, a 5-cycle $C_5$.

% -- independent sets
\begin{theorem}
\label{IE2} Given a $k$-chromatic graph $G$, $\{u,v\}\in V(G)$ is:
\begin{enumerate}
  \item an implicit-edge \textbf{iff} there is no independent set $S \in I(G-uv)$ such that $\{u,v\}\in S$ and $\chi(G-S)<k$.
  \item an implicit-identity \textbf{iff} there is no independent set $S \in I(G-u)$ such that $v \in S$ and $\chi(G-S)<k$.
\end{enumerate}
\begin{align}
      \text{(1)}& \; \left\{ S \in I(G-uv) \; | \; \{u,v\} \in S, \; \chi(G-S) < \chi(G) \right\} &= \emptyset \\
      \text{(2)}& \; \left\{ S \in I(G-u) \; | \; v \in S, \; \chi(G-S) < \chi(G) \right\} &= \emptyset
\end{align}
\begin{proof} $ $\\
\begin{enumerate}
  \item proof:
  \begin{enumerate}
  \item Assume that $\{u,v\}$ is an implicit-edge, but suppose that there is an independent set $S \in I(G-uv)$ where $\{u,v\}\in S$ such that $\chi(G-S)<k$. Then, we can find a ($k$-1)-coloring of $G-uv-S$ and restore all the vertices in $S$ with color $k$. Hence there will be a $k$-coloring of $G-uv$ where $c(u)=c(v)$, which is a contradiction.
  \item Assume that $\{u,v\}$ is \textbf{not} an implicit-edge. Then, by definition, there is a $k$-coloring of $G-uv$ such that $c(u)=c(v)$. Hence, there will be an independent set $S \in I(G-uv)$ where $\{u,v\}\in S$ such that $\chi(G-S)<k$.
  \end{enumerate}
  \item proof:
  \begin{enumerate}
  \item Assume that $\{u,v\}$ is an implicit-identity, but suppose that there is an independent set $S \in I(G-u)$ such that $v \in S$ and $\chi(G-S)<k$. Then, we can find a ($k$-1)-coloring of $G-S$ and restore all the vertices in $S$ with color $k$. Hence there will be a $k$-coloring of $G$ where $c(u)\neq(v)$ since $c(v)=k$, which is a contradiction.
  \item Assume that $\{u,v\}$ is \textbf{not} an implicit-identity. Then, by definition, there is a $k$-coloring of $G$ such that $c(u)\neq c(v)$. Hence, there will be an independent set $S \in I(G-u)$ such that $v \in S$ and $\chi(G-S)<k$.
\end{enumerate}
\end{enumerate}
\end{proof}
\end{theorem}

\begin{definition}
Given a $k$-chromatic graph $G$, we say that an independent set $S \subset V(G)$ is critical \textbf{iff} $\chi(G-S)= k-1$ (i.e. $S$ is a color class of some $k$-coloring of $G$)
\end{definition}
From this definition it follows immediately the that critical independent sets does not contain implicit-edges, due to theorem~\ref{IE2}. Hence every maximal independent set (also a minimal dominating set) containing an implicit-edge or containing exactly one vertex of an implicit-identity is not critical.

% -- invariant
\begin{theorem}
\label{ieinvariant}
If there is an implicit-relation (edge or identity) $\{u,v\}$ in a $k$-chromatic graph $G$ and $H=G-S$ is the resulting ($k$-1)-chromatic graph obtained by removing a critical independent set $S$ from $G$ then $\{u,v\}\in V(H)$ has the same implicit-relation (edge or identity, respectively) in $H$.

\begin{proof} $ $\\
\begin{enumerate}
  \item Let us suppose that $\{u,v\}$ is \textbf{not} an implicit-edge of $H$. Now, find a ($k$-1)-coloring of $H$ such that  $c(u)=c(v)$ and then restore the critical independent set $S$ with all vertices colored the same. Thus, the graph $G=H+S$ has a $k$-coloring where $c(u)=c(v)$ which is a contradiction.
  \item Let us suppose that $\{u,v\}$ is \textbf{not} an implicit-identity of $H$. Now, find a ($k$-1)-coloring of $H$ such that $c(u)\neq c(v)$ and then restore the critical independent set $S$ with all vertices colored the same. Thus, the graph $G=H+S$ has a $k$-coloring where $c(u)\neq(v)$ which is a contradiction.
\end{enumerate}
\end{proof}
\end{theorem}
%
%\begin{theorem}
%Let $H$ be a $k$-chromatic graph and $\{u,v\}\in V(G)$ implicit-edge of $G$
%Para que $G$+algo tenga al menos el mismo arco implicit entonces G+algo-S debe tener el mismo arco implicito.
%
%Lo reducimos al caso 2 chromatico, si tenemos un grafo B 2 chromatico y un arco implicito u,v  si añadimos mas vertices V2 y queremos que u,v siga siendo arco implicito entonces tenemos que conectar los vertices de forma tal que B+V2-S siga teniendo los caminos que unen u,v de orden par y sólo debe haber de caminos de orden par.
%\end{theorem}

% -- generalized invariant
\begin{theorem}
\label{IE3}
If there is an implicit-relation (edge or identity) $\{u,v\}$ in a $k$-chromatic graph $G$ then $\{u,v\}\in V(H)$ has the same implicit-relation (edge or identity, respectively) in $H$, for all subgraphs $H \subset G$ such that:
$$ H = G-\{S_1,\ldots,S_{\ell}\}\; \text { where } 0 \leq \ell \leq k-2, $$
where each $S_i$ a critical independent set not containing $u$ nor $v$.

\begin{proof} $ $\\
By reiterative application of theorem~\ref{ieinvariant} we get an inductive proof with base case $\ell$=0 and theorem~\ref{ieinvariant} for $\ell+1$.
\end{proof}

\end{theorem}

\begin{theorem}
\label{iekch}
If $G$ is a $k$-chromatic graph with an implicit-relation $\{u,v\}\in V(G)$:
\begin{enumerate}
  \item if $\{u,v\}$ is an implicit-edge then every $k$-coloring of $G$ contains at least one chain $\Omega_{uv}$ with colors $c(u),c(v)$ where $\{u,v\}\in \Omega_{uv}$
  \item if $\{u,v\}$ is an implicit-identity then every $k$-coloring of $G$ contains at least ($k$-1)-chains $\Omega_{ui}$ with colors $c(u),c(i)$  where $\{u,v\} \in \Omega_{ui} \; \forall i=1...(k-1)$.
\end{enumerate}

\begin{proof} $ $ \\
\begin{enumerate}
  \item Otherwise we can simply flip the color of vertex $u$ with other vertex $w$ of color $c(v)$ obtaining a $k$-coloring of $G-uv$ where $c(u)=c(v)$ which is a contradiction since $\{u,v\}$ is an implicit-edge.
  \item Otherwise we can simply flip the color of vertex $u$ with some
       adjacent vertex $i$ and, since there is no $\Omega_{ui}$ chain, there is a $k$-coloring of $G$ where $c(u)\neq c(v)$ but $\{u,v\}$ is an implicit-identity which is a contradiction.
\end{enumerate}
\end{proof}
\end{theorem}

\section{Chromatic Polynomials}

The chromatic polynomial $P(G,k)$ for a given graph $G$ is a
polynomial which encodes the number of different $k$-colorings of $G$, we can denote the number of $k$-colorings of a
graph $G$ as $P(G,k)$. Chromatic Polynomials satisfies the next two relations~\cite{Dong2005}(pp.4--6):

\begin{align}
\label{eq1} P(G,k)&= P(G-e,k)-P(G/e,k) \\
\label{eq2} P(G,k)&= P(G + e,k) +  P(G/e,k)
\end{align}

A direct consequence of (\ref{eq1}), (\ref{eq2}) and the definition of implicit-relations is presented in two new theorems:

\begin{theorem}\label{ie:tdcr} If $G$ is a $k$-chromatic graph and $\{u,v\}$ is an
implicit-edge of $G$ then $P(G/x,y)=0$, that is, the graph  $G/x,y$ is not $k$-colorable.
\begin{proof}$ $\\
If $e=\{u,v\}$ an implicit-edge of $G$ the it follows that $P(G,k) = P(G-e,k)$ or $P(G,k) = P(G+e,k)$, hence:
\begin{align}
P(G,k) &= P(G-e,k) \Longrightarrow  P(G/e,k)=0 \; \text{ due to (\ref{eq1})}. \\
P(G,k) &= P(G+e,k) \Longrightarrow  P(G/e,k)=0 \; \text{ due to (\ref{eq2})}.
\end{align}
\end{proof}
\begin{corollary} Given a minor-closed class $\;\mathcal{G}$ of graphs (i.e. if $G
\in \mathcal{G}$ and $G \geq H$, then $H\in \mathcal{G}$). If $G$ is a $k$-chromatic graph in $\mathcal{G}$ and $\{e\}$ a drawn-implicit-edge then $G/e$ is a ($k$+1)-chromatic graph in $\mathcal{G}$.
\end{corollary}
\end{theorem}

\begin{theorem} \label{ii:tdcr} If $G$ is a $k$-chromatic graph and $\{u,v\}$ is an
implicit-identity of $G$ then $P(G,k) = P((G-uv/u,v),k)$.
\begin{proof} $ $\\
Being $\{u,v\}$ an implicit-identity of a $k$-chromatic graph $G$ it follows that $P(G+e,k)=0$ hence:
\begin{equation}
P(G+e,k)=0 \Longrightarrow  P(G,k) = P((G-uv/u,v),k) \; \text{ due to (\ref{eq2})}.
\end{equation}
\end{proof}
\end{theorem}

\section{Planar graphs}
A graph is called planar if it can be drawn in a plane without edge crossings. Planar graphs poses the next important properties:
\begin{itemize}
    \item Every planar graph is four colorable~\cite{ah1,Robertson}.
    \item Planar graphs are closed under edge contraction, that is, every edge contraction of a planar graph results in a new planar graph.
\end{itemize}

\begin{figure}[h]
\begin{centering}
  \begin{tabular}{cc}
   \includegraphics[scale=1]{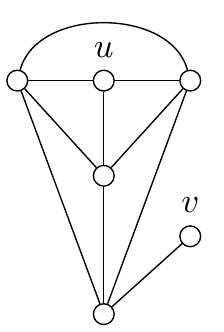} & \includegraphics[scale=1]{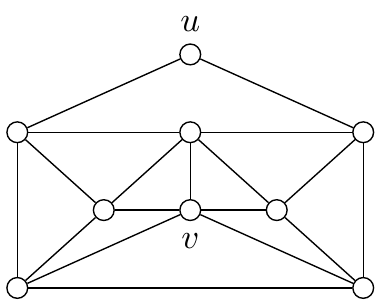} \\
   (A) & (B)
 \end{tabular}
  \caption{Implicit edges ($u,v$) in 4-chromatic planar graphs.}
  \label{planar4}
\end{centering}
\end{figure}

So, are there implicit-relations in 4-chromatic planar graphs? The question is affirmative, there are implicit-relations in
4-chromatic planar graphs, as can be seen in the graphs of Fig.~\ref{planar4}.

The main observation from Fig.~\ref{planar4}. is that the implicit-edge can not be drawn without edge line crossings. It is
straightforward to proof that if implicit-edges in 4-chromatic planar graphs could be drawn without edge line
crossings then this would be a counter example to the four-color theorem due to the corollary of Theorem~\ref{ie:tdcr}. Even more, the problem of determining if there are no implicit-edges in 4-chromatic planar graphs without edge-lines crossings is equivalent to a proof of the four-color theorem and conversely:

\begin{theorem} \label{ii:4cthm} There is a 5-chromatic planar graph \textbf{iff} exist a 4-chromatic planar graph $G$ such that $u,v\in V(G)$ is an implicit-edge and $G/u,v$ is planar.
\begin{proof}$ $\\
\begin{enumerate}
  \item Let $G$ be a critical 5-chromatic planar graph. Then by doing an edge-subdivision $H=G-e+w+uw+wv$ the resulting edges $uw$ and $wv$ of the 4-chromatic planar graph $H$ are implicit-edges and have no edge-line crossings.
  \item The backward implication is immediate from theorem~\ref{ie:tdcr}.
\end{enumerate}
\end{proof}
\end{theorem}
Finally, can it be proved, without using the 4-colors theorem, that:
\begin{theorem}\label{ie:noplanar} If $G$ is a 4-chromatic planar graph such that $u,v\in V(G)$ is an implicit-relation (edge or identity) then $G+uv$ is not a planar graph.
\end{theorem}

\section{Critical graphs}

\begin{theorem} If $G$ is a $k$-chromatic graph such that $w \in V(G)$ is critical then:
\begin{enumerate}
  \item if $\{u,v\}$ is an implicit-identity then $uw\in E(G)$ and $vw\in E(G)$.
  \item if $\{u,v\}$ is an implicit-edge then at least $uw\in E(G)$.  
\end{enumerate}
\begin{proof}$ $\\
\begin{enumerate}
  \item Otherwise $G-u-w$   is ($k$-1)-chromatic which also contradicts theorem~\ref{IE2}.
  \item Otherwise $G-u-v-w$ is ($k$-1)-chromatic which contradicts theorem~\ref{IE2}.  
\end{enumerate}
\end{proof}
\end{theorem}

\subsection{Double-critical graphs}
A $k$-chromatic graph $G$ is called double-critical if for every edge $uv \in E(G)$ the graph $G-u-v$ is $(k-2)$-colorable. A deep structural result about double-critical graphs can be proved easily with the help of implicit-relations:
\begin{theorem}
If $G$ is a $k$-chromatic double-critical graph and $uv\in E(G)$ then $u$ and $v$ have at least $k$-1 common neighbors, i.e.:
\begin{equation}
\chi(G-u-v)=\chi(G)-2 \rightarrow |N(u)\cap N(v)| \geq \chi(G)-2, \text{ for all } uv\in E(G)
 \end{equation}
\begin{proof}
Since in the ($k$-1)-chromatic graph $G-uv$, $\{u,v\}$ is an implicit-identity then, by theorem~\ref{iekch}, there will be $k-2$ color-chains containing $\{u,v\}$. Now since, $G$ is double-critical there is a ($k$-1)-coloring of $G-uv$ such that $\{u,v\}$ are the only vertices having the color $k$-1. Hence each color-chain can only pass trough a common neighbor, i.e. one neighbor for each color-chain.
\end{proof}
\end{theorem}